\documentclass[11pt,reqno,tbtags]{amsart}
\usepackage{amssymb}
\numberwithin{equation}{section}

\newtheorem{theorem}{Theorem}[section]
\newtheorem{lemma}[theorem]{Lemma}

\theoremstyle{definition}

\newtheorem{remark}[theorem]{Remark}

\newtheorem*{acks}{Acknowledgements}

\theoremstyle{remark}

\newenvironment{romenumerate}{\begin{enumerate}% gives (i), (ii) etc.
  }{\end{enumerate}}

\newcounter{thmenumerate}

\newcounter{xenumerate}   %no left indentation; thus wider lines

\newcommand{\refT}[1]{Theorem~\ref{#1}}

\newcommand{\refL}[1]{Lemma~\ref{#1}}
\newcommand{\refR}[1]{Remark~\ref{#1}}
\newcommand{\refS}[1]{Section~\ref{#1}}
\newcommand{\refand}[2]{\ref{#1} and~\ref{#2}}

\newcommand\marginal[1]{\marginpar{\raggedright\parindent=0pt\tiny #1}}

\begingroup
   \divide\count255 by 60
   \multiply\count255 by -60
   \advance\count255 by \time
   \ifnum \count255 < 10 \xdef\klockan{\the\count1.0\the\count255}
   \else\xdef\klockan{\the\count1.\the\count255}\fi
\endgroup

\newcommand\set[1]{\ensuremath{\{#1\}}}

\newcommand\xxpar[1]{\left(#1\right)}
\newcommand\bigpar[1]{\bigl(#1\bigr)}
\newcommand\Bigpar[1]{\Bigl(#1\Bigr)}

\def\rompar(#1){\textup(#1\textup)}    % usage: \rompar(...)

\newcommand\parfrac[2]{\Bigpar{\frac{#1}{#2}}}

\def\xexp(#1){e^{#1}}

\newcommand\ntoo{\ensuremath{{n\to\infty}}}

\newcommand\ie{i.e.\spacefactor=1000}
\newcommand\eg{e.g.\spacefactor=1000}

\newcommand\whp{\textbf{whp}}

\newcommand\E{\operatorname{\mathbb E{}}}
\renewcommand\P{\operatorname{\mathbb P{}}}
\newcommand\Var{\operatorname{Var}}

\newcommand\Tr{\operatorname{Tr}}

\newcommand\ga{\alpha}
\newcommand\gb{\beta}
\newcommand\gd{\delta}

\newcommand\gam{\gamma}

\newcommand\gl{\lambda}

\newcommand\cE{\mathcal E}
\newcommand\cG{\mathcal G}
\newcommand\cH{\mathcal H}

\newcommand\ett[1]{\boldsymbol1[#1]}

\def\[#1]{[\![#1]\!]}

\newcommand\smatrixx[1]{\left(\begin{smallmatrix}#1\end{smallmatrix}\right)}
\newcommand\matrixx[1]{\left(\begin{matrix}#1\end{matrix}\right)}

\newcommand\qq{^{1/2}}

\renewcommand{\=}{:=}

\newcommand\lhs{left hand side}

\newcommand\ij{_{ij}}
\newcommand\io{_{i0}}
\newcommand\Gx{G^*}
\newcommand\Bx{B^*}
\newcommand\tB{\overline B}
\newcommand\gcnq{G_c(n,2q)}
\newcommand\gcnqx{\Gx_c(n,2q)}
\newcommand\gcnqiix{\Gx_c(n,2q-2)}
\newcommand\bnqx{\Bx(n,n;2q)}

\newcommand\REM[1]{\texttt{[#1]}\marginal{XXX}}

\newcommand\urladdrx[1]{\urladdr{\def~{\~{}}#1}}

\def\k{\kappa}
\begin{document}
\title%[]
{Rainbow Hamilton cycles in random regular graphs}

\date{August 1, 2005} %October 19, 2004}

\author{Svante Janson}
\address{Department of Mathematics, Uppsala University, PO Box 480,
SE-751 06 Uppsala, Sweden}
\email{svante.janson@math.uu.se}
\urladdr{http://www.math.uu.se/~svante/}

\author{Nicholas Wormald}\thanks{The second author acknowledges the support of
the  Canadian Research Chairs Program and NSERC}
\address{Department of Combinatorics and Optimization,
University of Waterloo,
Waterloo ON,
Canada N2L 3G1}
\email{nwormald@uwaterloo.ca}
\urladdrx{http://www.math.uwaterloo.ca/~nwormald/}

%\keywords{<keywords>}
%{60C05 (68P10,68W40)} %%{Primary: <subject>; Secondary: <subject>}

\dedicatory{To Alan Frieze on the occasion of his 60th birthday}

\begin{abstract} 
A rainbow subgraph of an edge-coloured graph has all edges of distinct
colours.  A 
random $d$-regular graph with $d$ even, and having edges coloured
randomly with $d/2$ 
of each of $n$ colours, has a rainbow Hamilton cycle with probability
tending to 1 as 
$n\to\infty$, provided $d\ge 8$.
\end{abstract}

\maketitle

\section{Introduction}\label{S:intro}
An  edge-coloured graph is a
\emph{rainbow} if no colour appears more than once. We will study
rainbow Hamilton 
cycles in edge-coloured graphs with
$n$ vertices where the number of colours available is also $n$; thus a rainbow
Hamilton cycle uses each of the 
colours exactly once.
 
We consider in this paper only (random) regular graphs and
(random) colourings where each colour occurs the
same number of times. If each colour occurs on
$q$ edges, we thus have $qn$ edges, and hence the vertex degrees are
$2q$.   
We use the standard notation $G(n,d)$ for a uniformly chosen random
$d$-regular graph on $n$ given (labelled) vertices.
We will only consider the case $d=2q$ even. 
(Hence there is no parity restriction on $n$.) 
Having sampled a random graph $G(n,2q)$, we then randomly colour its
$qn$ edges by $n$ colours ($1,\dots,n$, say) with $q$ edges of each colour,
again choosing uniformly among all possibilities.
We denote the resulting randomly coloured random graph by $\gcnq$.

Our main result is the following on randomly coloured random regular
graphs.
(For some related results on the random graph $G(n,m)$, 
see Cooper and Frieze \cite{CF}.)
We say that an event holds \emph{with high probability} (\whp), 
if it holds with probability tending to 1 as
$n\to\infty$.
(All unspecified limits in this paper are for \ntoo.)

\begin{theorem}\label{T1}
Consider the randomly coloured random $2q$-regular graph $\gcnq$, 
with $n$ colours and
$q$ edges of each   colour. 
Then, \whp, there exists a rainbow Hamilton cycle if $q\ge4$,
  and  not if $q\le3$.
\end{theorem}

Recall that it was shown by Robinson and Wormald \cite{RW92,RW94} that
$G(n,d)$ \whp{} contains a Hamilton cycle as soon as $d\ge3$.
In our setting, when $d=2q$ has to be even, we thus \whp{}
have Hamilton cycles, ignoring the colouring, when $d\ge4$, but
rainbow Hamilton cycles only when $d\ge8$.  
It is nevertheless remarkable that \whp\ some Hamilton cycle manages
to pick up an edge of each colour in a random 8-regular graph, when there
are only four edges of each colour to choose from. 
 
\begin{remark}\label{Rdecom}
  In a similar direction, Robinson and Wormald~\cite{RW01}  showed
that a random 3-regular graph with $o(\sqrt n)$ randomly specified edges \whp\
has a Hamilton cycle passing through all the specified edges (and, moreover, in
randomly prespecified directions). It has further been shown by  Kim and
Wormald~\cite{KW99} that a random $2q$-regular graph \whp{} has an
edge-decomposition into $q$ Hamilton cycles, provided $q\ge2$.

It is natural to ask whether, similarly, a randomly coloured
$2q$-regular graph with $n$ colours and $q$ edges of each colour,
as above,
\whp{} has an
edge-decomposition into $q$ rainbow Hamilton cycles.
By computing the expected number of such
decompositions (similarly to the proof of \refL{LE} below),
it is easily seen
that this is \whp{} false when $q\le 4$.
We leave the case $q\ge5$ (when the expected number tends to infinity)
as an open problem.
\end{remark}

The proof of \refT{T1} is based on the small subgraph
conditioning method introduced
by Robinson and Wormald \cite{RW92,RW94}, and further developed in
\cite{SJ103}, \cite{MRRW}, \cite{NW99} and
\cite[Chapter 9]{JLR}.
However, for this problem we have to consider the colourings of the
small subgraphs too, see \refS{S:small}.

\begin{acks}
   This problem was suggested by Alan Frieze during the Conference on
   Random Structures and Algorithms at Emory University, Atlanta, 1995;
he was originally intended as a coauthor but later declined this.
Consequently we are pleased to be able finally to dedicate this work to him,
marking 10 years since its beginnings, in which he was involved, and 60
since his.

We also acknowledge the assistance of the Maple algebraic manipulation
package for the variance calculations in Section~\ref{s:var}.
Although the proof we found can be verified by hand, Maple was
instrumental in finding that proof.  
\end{acks}

\section{Multigraphs, a bipartite graph, and traffic rules}
\label{S:prel}

As usual in the study of random regular graphs (with small degree), it
is convenient to extend the study to multigraphs.
Recall that a convenient way (at least for theoretical purposes) 
to generate a random regular graph is the so-called
\emph{configuration model} or \emph{pairing model},
see \eg{}~\cite{bollobas} or~\cite{NW99}: 
We start with $nd$ points partitioned into $n$ cells of $d$ points
each. We then take a random pairing of the points into $nd/2$ pairs
(assuming $nd$ to be even). Collapsing each cell to a vertex and
regarding each pair as an edge, we obtain a random $d$-regular
multigraph that may contain loops and multiple edges; we denote this
random multigraph by $\Gx(n,d)$. 
(The points themselves are called \emph{half-edges}.) 
It is well-known, and easily seen,
that if we condition $\Gx(n,d)$ on being a simple graph (no loops nor
multiple edges), then we obtain the uniformly distributed random
regular graph $G(n,d)$. Moreover, it is well-known
that for fixed $d$, the probability $\P(\Gx(n,d)\text{ is simple})$
tends to a non-zero limit as \ntoo; hence, every property that
$\Gx(n,d)$ has \whp, is  \whp\ enjoyed by $G(n,d)$ too.  
In particular, we may (and will) prove \refT{T1} by proving the
following extension of it; we define $\gcnqx$  by analogy with $\gcnq$,
by choosing uniformly at random
a colouring of the edges with   $n$ colours with $q$ edges of each colour.

\begin{theorem}\label{T2}
\refT{T1} holds for the randomly coloured random regular multigraph
$\gcnqx$
too.
\end{theorem}

We find it useful to introduce an associated   bipartite graph.
(This is really a multigraph too, since it may have multiple edges.)
Given the randomly coloured multigraph $\gcnqx$, add a new vertex on
each edge. We give each new vertex the colour of the edge it bisects,
leaving the original vertices uncoloured.
Finally, we combine the $q$ vertices of each colour into a single coloured
vertex of that colour. This gives us a $2q$-regular
bipartite (multi)graph with $n+n$ vertices; the original $n$ vertices form
one side of the bipartition, and the $n$ coloured vertices the other.
Moreover, each coloured vertex comes with a pairing of the edges (or
half-edges) attached to it; this pairing shows which pairs of edges
  correspond to edges in the multigraph.
We may think of the coloured vertices as having $2q$ attached
half-edges arranged in a circle, with each half-edge matched  to the
opposite one. 
There is then a one-to-one correspondence between walks in the
multigraph and walks in the bipartite graph
(of twice the length, and beginning and ending at blank   vertices)
that pass `straight ahead' between matched half-edges at  each
coloured vertex. In particular, rainbow Hamilton 
cycles in the multigraph correspond to Hamilton cycles in the bipartite graph
that obey this
 traffic rule.

The colours are no longer important in the bipartite graph, but it
will be convenient to refer to the two sets of vertices in the
bipartition as `plain' and `coloured'.

Conversely, we may start with the bipartite graph, with given traffic
rules, and obtain the original multigraph by combining the edges two
by two at the coloured vertices. Note that choosing the bipartite
(multi)graph at random using the configuration model (in its bipartite
version, and with traffic rules as above given in each coloured cell)
gives back 
the random coloured multigraph $\gcnqx$ with the right distribution.
% We let $\bnqx$ denote this random bipartite multigraph with traffic
% rules as above.
% To prove Theorems \refand{T1}{T2}, it is thus enough to prove the following.
% \begin{theorem}\label{T3}
%Consider the random bipartite multigraph $\bnqx$ with $n+n$
%vertices of degree $2q$, and with a pairing of the half-edges at each
%vertex in one   part of the bipartition.
%  Then, \whp, there exists a rainbow Hamilton cycle obeying the
 
We let $\bnqx$ denote this random bipartite multigraph with traffic
rules as above,   that is,  with $n+n$
vertices of degree $2q$ and with a random pairing of the half-edges
at each vertex in the coloured part of the bipartition.   (Although
viewing it as a multigraph when referring to cycles etc., all
computations are done with the equivalent configuration model.) 
To prove Theorems \refand{T1}{T2}, it is  enough to prove the following.
 
\begin{theorem}\label{T3}
The random bipartite multigraph $\bnqx$ \whp\ has a rainbow Hamilton
cycle obeying the traffic rules  if
$q\ge4$, and not if $q\le3$.
\end{theorem}

\begin{remark}
One may study random (regular) graphs with other 
traffic rules at the vertices.
In general, we may equip each vertex of degree $d$ with a 
(possibly directed, and possibly random)
\emph{connection graph} 
with $d$ vertices representing the incident edges; the edges 
in the connection graph show the allowed connections between incoming
and outgoing edges.
In our case, the connection graph is the complete graph (no
restrictions) for one side of the bipartition, and a matching with $d/2$
edges for the other side.
We do not know of any general study, but a few examples of
this type have appeared in the literature:

Garmo \cite{Garmo95,Garmo99} studied random railways; these
are regular (typically cubic) graphs where the vertices (representing
switches) have 
connection graphs that are stars.
%Some extensions to other connection graphs have been made by Pontus
%Andersson (personal communication).
In \cite{GarmoJK}, this was extended to graphs where a random subset
of the vertices have a star as connection graph and the rest the
complete graph. 

Gamburd \cite{Gamburd} studied long cycles in random oriented
cubic graphs; here the connection graph is a directed 3-cycle at each vertex.
\end{remark}

\section{Small subgraphs}\label{S:small}

The small subgraph conditioning method introduced
by Robinson and Wormald \cite{RW92,RW94}
has been successfully applied
to several problems, in particular in the theory of random regular
graphs, see \eg{} \cite[Chapter 9]{JLR}, \cite{NW99} and \cite{GJKW}.
(For applications to random hypergraphs, see \cite{FJ,CFMR}.)

As often pointed out by Alan Frieze, see \cite{CFMR,FJ,FJMRW},
the method can be regarded as an \emph{analysis of variance}.
The main idea is that we consider some random variable, $Y$ say,
that counts occurrences of some
structure, and let a parameter
\ntoo.  Typically, it is easy to
prove that the expectation $\E Y$ tends to infinity, but we want to
show that $\P(Y>0)\to1$.
If the variance $\Var(Y)$ is $o(\E Y)^2$, then the second moment
method (\ie{} Chebyshev's inequality) immediately shows the desired result.
The small subgraph conditioning method applies to cases where
the variance $\Var(Y)$ is of the same order as $(\E Y)^2$,
by showing that the variance can be explained, up to a factor $1-o(1)$,
by the interaction between the
numbers of some small subgraphs and the random variable $Y$. The desired
conclusion
$Y>0$
\whp{} then follows by conditioning on the numbers of these small
subgraphs and using Chebyshev's inequality on the conditioned variables.
For details, see \cite[Theorem 9.12--Remark 9.18]{JLR} and
\cite[Theorem 4.1]{NW99}. We state the results there in the following
form (an immediate consequence of \cite[Corollary 4.2]{NW99}).
We use $[x]_m\=x(x-1)\cdots(x-m+1)$ to denote falling factorials.

\begin{theorem}\label{T0}
Let $\lambda_i>0$ and $\delta_i\geq -1$ be real numbers for $i=1,2,\ldots$
and suppose that for each $n$ there are random variables $X_i =
X_i(n)$, $i=1,2,\ldots$
and $Y=Y(n)$, all defined on the same probability space $\cG = \cG_n$
such that $X_i$ is nonnegative integer valued, $Y$ is nonnegative and
$\E Y>0$ (for $n$ sufficiently large).  Suppose furthermore that
\begin{romenumerate}
\item For each\, $k\geq 1$, the variables $X_1,\ldots ,X_k$ are asymptotically
independent Poisson random variables with $\E X_i\rightarrow\lambda_i$,
\item if $\mu_i=\lambda_i(1+\delta_i)$, then
   \begin{equation}\label{sofie}
\frac{\E(Y[X_1]_{m_1}\cdots [X_k]_{m_k})}{\E Y}
  \rightarrow \prod_{i=1}^k \mu_i^{m_i}
   \end{equation}
for every finite sequence $m_1,\ldots ,m_k$ of nonnegative
integers,
\item $\sum_i \,\lambda_i\,{\delta_i}^2 < \infty$,
\item $\E Y^2/(\E Y)^2\leq \exp(\sum_i \lambda_i\,{\delta_i}^2) + o(1)$
as $n\rightarrow\infty$.
\end{romenumerate}
Then, if $\cE$ is the event $\wedge_{\delta_i=-1} \set{X_i=0}$,
$\P(Y>0\mid\cE)\to1$.
In particular, if $\gd>-1$ for every $i$, then $Y>0$ \whp.
\end{theorem}

We will actually use \refT{T0} with a doubly indexed sequence $X\ij$;
obviously, this is just a matter of notation.

In many applications of the small subgraph conditioning method, the
variables $X_i$ are the numbers of cycles of different lengths. This
has   perhaps misled some into the belief that the short cycles are
expected to provide the answer in all cases. But they play the central
role for most problems only because they are the only possible
`unusual' small subgraphs. The subgraphs of fixed size in a random
$d$-regular graph are very well behaved. Near a random vertex, such a
graph looks locally like a tree. But even that statement can be
misleading when we consider what comes shortly. The thing to focus on
is that, because \whp\ no two short cycles are near each other, the
number of subgraphs of any particular type are determined by the
numbers of short cycles. 

In our case, it will turn   out that conditioning on the
numbers of small cycles in
$\gcnqx$ does not explain all of the variance of the number of rainbow Hamilton
cycles; we have to consider also colourings. Note that for each fixed length
$i$, there are only a few cycles of length $i$ (the expected number is $O(1)$),
and \whp{} they are all rainbow, so   we would expect no explanation of
variance to be caused by the numbers of intrinsically differently coloured
short cycles. However, we may consider, for example, the number of paths of
length
$i$ where the first and last edges have the same colour. The expected number is
$\Theta(1)$.  These are intrinsically different from short rainbow
paths, and it turns out that these structures too will be significant in the
analysis of variance.  

Perhaps surprisingly, there is even more to consider. The existence of
two short 
paths, each joining a blue edge to a red edge, is significant, even though they
are two different blue edges and two different red edges in distant
parts of the 
graph. However, this is not surprising given the discussion above
about small subgraphs. Colours are clearly relevant in our present
problem, so we should consider coloured subgraphs. Typical small
subgraphs, not necessarily neighbourhoods of vertices, are forests
with distinctly coloured edges. The numbers of small forests in which
some of the edges are coloured the same thus qualify as special small
subgraphs. It was for this reason that the method was called small
subgraph conditioning  in~\cite{NW99}, rather than short cycle
conditioning. This is indeed the first application of the method in
which the small subgraphs involved are disconnected.   

 To describe the general situation precisely, we work with the random
bipartite multigraph
$\bnqx$ defined in \refS{S:prel}, and let $Y$ be the number of Hamilton cycles
 in the multigraph that obey the traffic rules. 
Recall that $Y$ equals the number of rainbow Hamilton cycles in $\gcnqx$.
Further, for each $i\ge1$ and $j$ with $0\le j\le i$,
we let $X\ij$ be the number of cycles of length $2i$ in
$\bnqx$ that violate the traffic rules at exactly coloured $j$ vertices.
(Thus, $Y=X_{n0}$, but we are mainly interested in $X\ij$ for small
$i$.)
Note that $X_{i0}$ equals the number of rainbow $i$-cycles in $\gcnqx$, 
and thus it \whp{} equals the number of $i$-cycles in $\gcnqx$,
while $X_{i1}$ \whp{} equals the number of paths of length $i+1$ where
the first and last edges have the same colour. (This holds only \whp,
since for $X_{i1}$ the endpoints of the path may coincide with each
other or with some interior point.) We may similarly interpret $X\ij$
for $j\ge2$, at least \whp, as the number of certain collections of
$j$ paths,   generalising the example mentioned above, but we leave
the details to the reader.

We state three lemmas that will be proven in the following sections.

 \begin{lemma}
  \label{LE}
Suppose that $d=2q\ge4$. Then
  \begin{equation*}
\E(Y) 
=
\Theta\xxpar{\frac{(d-1)(d-2)^{d-2}}{d^{d-2}}}^n.
  \end{equation*}
Hence, as \ntoo, $\E(Y)\to0$ for $d\le 6$ but
$\E(Y)\to\infty$ for $d\ge 8$.
\end{lemma}

\begin{lemma}
  \label{L1}
Conditions \textup{(i)} and 
\textup{(ii)} in \refT{T0} are satisfied for
the variables 
$(X\ij)_{ij}$ and 
\begin{align*}
 \gl\ij&= 
\frac{1}{2i}\binom{i}{j}(d-1)^i(d-2)^j,
\\
\gd\ij&=
\begin{cases} 
(-1)^{i+j}\frac{2^j}{(d-1)^i(d-2)^j}, & j>0,
\\
-\frac{2}{(d-1)^i}\ett{i \text{ odd}}, & 
j=0.
\end{cases}
\end{align*}
\end{lemma}

  \begin{lemma}
   \label{LE2}
Suppose that $d>4$. Then
   \begin{equation*}
\E Y^2/(\E Y)^2 \to	\parfrac{d}{d-4}^{1/2}.
   \end{equation*}
\end{lemma}

\begin{proof}[Proof of \refT{T3}]
First note that if $d\le6$, then $\E(Y)\to0$ by \refL{LE}, and thus 
$\P(Y>0)\to0$, \ie{} $Y=0$ \whp. In other words, there is then \whp{} no
rainbow Hamilton cycle in $\gcnqx$.

For the remainder of the proof, assume that $d=2q\ge8$. By \refL{LE},
$\E(Y)\to\infty$. We want to show that $Y>0$ \whp.
We employ \refT{T0} with $X\ij$ as defined above, and $\gl\ij$ and
$\gd\ij$ as given in \refL{L1}. Note that $\gd\ij>-1$ for all $i$ and $j$,
so it remains only to show that the assumptions (i)--(iv) in \refT{T0}
hold.
For (i) and (ii), this is \refL{L1}.

For (iii) and (iv) we split the sum into two parts.
\begin{align*}
  \sum_{i,j>0} \gl\ij\gd\ij^2
&=
\sum_{i,j\ge1} \frac1{2i}\binom{i}j\frac{4^j}{(d-1)^i(d-2)^j}
\\&
=
\sum_{i=1}^\infty 
 \frac1{2i}(d-1)^{-i}\Bigpar{\Bigpar{1+\frac4{d-2}}^i-1}
\\&
=\frac12
\sum_{i=1}^\infty 
 \frac1{i}
\Bigpar{\parfrac{d+2}{(d-1)(d-2)}^i-\frac1{(d-1)^{i}}}
\\&
=-\frac12\ln\Bigpar{1-\frac{d+2}{(d-1)(d-2)}}
 +\frac12\ln\Bigpar{1-\frac{1}{d-1}}
\\&
=-\frac12\ln\frac{d^2-4d}{(d-1)(d-2)}
 +\frac12\ln\frac{d-2}{d-1}
\\&
=\frac12\ln\frac{(d-2)^2}{d(d-4)}
\end{align*}
and
\begin{align*}
  \sum_{i} \gl_{i0}\gd_{i0}^2
&=
\sum_{i\text{ odd}} \frac1{2i}\frac{4}{(d-1)^i}
=-\ln\Bigpar{1-\frac{1}{d-1}}
 +\ln\Bigpar{1+\frac{1}{d-1}}
\\&
= -\ln\frac{d-2}{d-1} 
 +\ln\frac{d}{d-1}
=\ln\frac{d}{d-2}
\end{align*}
Consequently,
\begin{align*}
 \sum_{i,j} \gl\ij\gd\ij^2
=\frac12\ln\frac{(d-2)^2}{d(d-4)}
+\ln\frac{d}{d-2}
=\frac12\ln\frac{d}{d-4}.
\end{align*}
This proves (iii), and together with \refL{LE2} also (iv).
\end{proof}

\section{Expectation}\label{SE}

\begin{proof}[Proof of \refL{LE}]
There are $n!^2/2n$ ways to arrange the $2n$ vertices in a cycle, 
with plain and coloured vertices alternating,
and for each such arrangement $d(d-1)$ ways to choose the half-edges at
each plain vertex and $d$ ways to choose the half-edges at each
coloured vertex (obeying the traffic rules). For each such choice, the
probability that the selected $4n$ half-edges are connected to each other
in the specified order equals $((d-2)n)!/(dn)!$.
Consequently,
using Stirling's formula,
\begin{equation*}
  \E Y = \frac{d^{2n}(d-1)^n n!^2 ((d-2)n)!}{2n\,(dn)!}
=\Theta\Bigpar{\frac{d^2(d-1)(d-2)^{d-2}}{d^{d}}}^n
=\Theta\bigpar{f(d)}^n,
\end{equation*}
where $f(d)\=(d-1)(1-2/d)^{d-2}$. 
We have $f(4)=3/4<1$,
$f(6)=80/81<1$,
$f(8)=5103/4096>1$, and $f(d)>(d-1)e^{-2}>1$ for $d>8$.
\end{proof}

\section{Short  cycles} 
\label{SC}
\begin{proof}[Proof of \refL{L1}]
We use arguments that have become standard for similar problems for
random ragular graphs, see \eg{} \cite[Section 9.4]{JLR} or
\cite[Section 4.2]{NW99}; we will thus omit some details.
	
For (i), we use the method of moments. It suffices to show that
\begin{equation*}
\E \prod \ij [X\ij]_{m\ij}\to\prod_{i=1}^k\gl\ij^{m\ij}  
\end{equation*}
for every finite set of non-negative integers \set{m\ij}.
For convenience, we will only treat the expectation of a single $X\ij$;
as in all similar problems, the argument extends immediately to
(mixed) higher factorial moments.

To calculate $\E X\ij$, we count   the appropriate oriented cycles
with a designated initial vertex, which we require to be plain; this
counts each 
cycle
$2i$ times. The vertices in the cycle may now be chosen in 
$[n]_i^2\sim n^{2i}$ ways. 
Consider first the case $j=0$, \ie{} cycles obeying the traffic rules
everywhere. 
For each choice of vertices there are, as in \refS{SE},
$d(d-1)$ ways to choose the half-edges at
each of the $i$ plain vertices and $d$ ways to choose the half-edges at each
of the $i$ coloured vertex. Finally, the probability of pairing the $4i$
chosen half-edges into $2i$ edges is $1/[dn]_{2i}\sim (dn)^{-2i}$.
Hence,
\begin{align*}
  2i \E X_{i0}
= \frac{d^{2i}(d-1)^i[n]_i^2}{[dn]_{2i}}
\to
(d-1)^i.
\end{align*}
For $j>0$ we argue similarly. The traffic rules are to be violated
at  precisely $j$ coloured vertices. These may be chosen in $\binom
ij$ ways, and at each of them there is additional factor of $d-2$ for
the choice 
of the out-going half-edge.
Hence we obtain, for all $i$ and $j$,
\begin{align*}
  2i \E X_{ij}
= \binom ij\frac{d^{2i}(d-1)^i(d-2)^j[n]_i^2}{[dn]_{2i}}
\to
\binom ij (d-1)^i(d-2)^j,
\end{align*}
or $\E X\ij\to\gl\ij$.

For (ii), we first observe that the left hand side of \eqref{sofie},
by symmetry, remains the same if we fix two half-edges at each vertex,
always choosing two opposite half-edges at the coloured vertices,
and then replace $Y$ by the indicator that the chosen half-edges
comprise a rainbow Hamilton cycle. Denoting this event by $\cH_1$, we
thus want to show 
\begin{equation}\label{erika}
\E \Bigpar{\prod \ij [X\ij]_{m\ij}\,\Big|\, \cH_1}
\to\prod_{i=1}^k\mu\ij^{m\ij}.
\end{equation}

Condition on $\cH_1$, and let $H_1$ be the (unique)
rainbow Hamilton cycle that uses the chosen half-edges. It is easily
seen that the remainder of the graph $\gcnqx$
can be regarded as the random multigraph 
$\gcnqiix$, and that this is independent of $H_1$.  
Hence, the \lhs{} of \eqref{erika} equals the
expectation in the union of a random rainbow Hamilton cycle $H_1$ and
an independent $\gcnqiix$ on the same vertex set.

For the same reasons as in (i),   we will only consider a single expectation
$\E(X\ij\mid\cH_1)$.
Consider first the case $j=0$. 
We may, as for (i), choose the vertices of the $2i$-cycle 
in   $n^{2i}(1+o(1))$ ways.
We then decide whether the $2i$ edges are in the Hamilton cycle $H_1$
or in $\gcnqiix$; we denote the choices by $\ga_s\in\set{1,2}$ for
$s=1,\dots,2i$. At a plain vertex where the incoming edge 
is to have type $\ga\in\set{1,2}$ and the outgoing edge type $\gb$,
there is for each possible incoming half-edge $a_{\ga\gb}$ choices of
the outgoing, where the numbers $a_{\ga\gb}$ are conveniently
collected in the matrix  
\begin{equation*}
 A\=(a_{\ga\gb}) = \matrixx{ 1 & d-2 \\ 2 & d-3}.
\end{equation*}
At the coloured vertices there is only one choice for the outgoing
edge, and it has to have the same type as the incoming; we encode this
as $b_{\ga\gb}$ with  
$ B\=(b_{\ga\gb}) = \smatrixx{ 1 & 0 \\ 0 & 1}
$, the identity matrix.
Note that we have not counted the number of incoming half-edges; this
is because these numbers cancel when we take into account the
probability of making the connections; the probability of connecting a
half-edge of either type to \emph{some} incoming half-edge of the same type
at a given vertex (in the opposite part of the bipartition)
is   $(1+o(1))n^{-1}$. Moreover, the probability that all $2i$
connections are made 
is $(1+o(1))n^{-2i}$, except in the case when all $\ga_s=1$, which is
impossible for $n>i$.
Consequently,
\begin{equation*}
  \begin{split}
 2i \E(X_{i0}\mid \cH_1)
&
\to \sum_{\ga_1,\dots,\ga_{2s}=1}^2 a_{\ga_1\ga_2} b_{\ga_2\ga_3} \dotsm
b_{\ga_{2i}\ga_1} -a_{11}^ib_{11}^i
\\&
=\Tr(AB)^i-1
=\Tr(A^i)-1	
=(d-1)^i+(-1)^i-1,
  \end{split}
\end{equation*}
since $A$ has the eigenvalues $d-1$ and $-1$.
We have thus shown \eqref{erika} for this case, with 
\begin{equation*}
\mu\io =\frac{(d-1)^i+(-1)^i-1}{2i}=\gl\io(1+\gd\io)
\end{equation*}
as required.

For $j>0$ we argue similarly. Now we have to choose $j$ coloured
vertices where the traffic rule is violated, and for these vertices
the matrix $B$ is replaced by
\begin{equation*}
 \tB\= A-B=\matrixx{ 0 & d-2 \\ 2 & d-4}.
\end{equation*}
Luckily, $A$ and $\tB=A-I$ commute, so all $\binom ij$ choices of the
violating vertices give the same result, and thus, for $j>0$,
\begin{equation*}
  \begin{split}
 2i \E(X\ij\mid \cH_1)
&
\to 
\binom ij \Tr(A^i\tB^j)
=\binom ij \Tr\bigpar{A^i(A-I)^j}	
\\&
=\binom ij \Bigpar{(d-1)^i(d-2)^j+(-1)^i(-2)^j},
  \end{split}
\end{equation*}
which equals $2i\gl\ij(1+\gd\ij)$ as required in this case too.
\end{proof}

\section{Variance} \label{s:var}

In this section we prove \refL{LE2}, thus completing the proof of \refT{T1}.
 
We   compute $\E Y^2$ by calculating the probability that a given  ordered
pair of Hamilton cycles $(H_1,H_2)$ are contained in the pairing
corresponding to $\bnqx$, and summing over all possible ordered $(H_1,H_2)$. 

 The first part is similar to the treatment of $H_1$ in Section~\ref{SC}. 
Let $k$ denote the number of coloured vertices in which the same half-edges
are used by both $H_1$ and $H_2$, and let $j$ denote the number of
blank  vertices of this type. If $H_1\ne H_2$ then the half-edges shared
by the two cycles occur in $k-j$   ``strings" of consecutive half-edges
around $H_1$. 
(Note that each string ends at plain vertices.)
The strings also occur in $H_2$, though in  a different order.
Together, $H_1$ and $H_2$ determine the pairs containing two of the
half-edges at $k$ coloured vertices and four at all other coloured vertices,
so $4n-2k$ pairs in all  (as each pair contains just one coloured vertex). 
Hence, defining $\cH_i$ as the event that $H_i$ occurs,
 $$
 \P( \cH_2\mid \cH_1) =  \frac{\big((d-4)n+2k\big)!}{\big((d-2)n\big)!},
$$
and we may write 
$$
\frac{\E Y^2}{(\E Y)^2}=\frac{1+\sum_{k,j<n}N(k,j)  \P ( \cH_2\mid \cH_1)}
{\E Y}
$$
where $N(k,j)$ is the number of different $H_2$ (or, to be precise, sets of
pairs corrsponding to $H_2$) overlapping any  given cycle $H_1$ with
particular values of
$k$ and
$j$. The term 1 accounts for the case $H_2=H_1$, when $k=j=n$.

Since $\E Y$ was given precisely in Section~\ref{SE}, all that is  left to
evaluate is $N(k,j)$. Note that the cardinality of the set $S_1$ of plain
vertices with three half-edges contained in $H_1\cup H_2$ is $2k-2j$, and
for the set $S_2$ of plain vertices with four half-edges, it is $n-2k+j$. 

By elementary counting, the number of ways to place $k-j$ strings  referred
to above onto $H_1$ is 
$$
\frac{n}{k}\binom{k}{j}\binom{n-k-1}{ k-j-1}.
$$
Here the first factor converts the problem from a cyclical one to a
linear one in which the first coloured vertex is one of the $k$
special  ones. The second factor is for 
deciding the relative positions of the $k$ coloured vertices in a
sequence of  $k-j$  strings and 
the third is for deciding the positions of the vertices not in strings. 

The
number of ways to choose the half-edges being used by $H_2$ at each plain
vertex is 
$ 
(d-2)^{|S_1|}((d-2)(d-3))^{|S_2|},
$
and $(d-2)^{n-k}$ for the half-edges at the set $S_3$ of coloured vertices
not used by
$H_1$  (obeying the traffic rules). These choices
determine the direction that $H_2$ passes through each vertex except for
those in $S_1$, which are determined by the direction it passes through each
string. These directions can be chosen in
$2^{k-j}$ ways, giving
$$
2^{k-j}(d-2)^{2n-k-j}(d-3)^{n-2k+j} 
$$
ways to make these choices. (For convenience we will count oriented versions
of $H_2$ and will divide by 2 at the end.)

Now that the order of half-edges used by $H_2$ is determined at each
vertex, it remains to choose the remaining $2(n-k)$ pairs. These
pairs must connect the strings and the vertices in $S_2
\cup S_3$ into a Hamilton cycle in a bipartite fashion (regarding each
string as a vertex) and connecting the `out' half-edge at a vertex to the
`in' one at the next. The number of such choices of pairs is
$$
 (n-k)!(n-k-1)!.
$$
Multiplying all the displayed factors together and dividing by 2 to un-orient
$H_2$ gives $N(k,j)$. Combining with the earlier equations then produces
\begin{equation}\label{mainvar}
\frac{\E Y^2}{(\E Y)^2} =\frac{1}{\E Y} +\sum_{k<n}\sum_{j<k}
f(n,d,k,j)
\end{equation}
where 
\def\Nkjnums{  n^2(k-1)!\,(n-k)!\,^3 2^{k-j}(d-2)^{2n-k-j}(d-3)^{n-2k+j}
  }
\def\Nkjdens{(n-k)^2 (k-j)!\,(k-j-1)!\,j!\,(n-2k+j)!}
\def\probnums{\bigl( (d -4)n+2k\bigr)!}
\def\probdens{\bigl((d-2)n\bigr)!\,^2}
\def\EYdens{  (dn)!}
\def\EYnums{d^{2n}(d-1)^n n!^2 }
\multlinegap=0pt
\begin{multline*} 
f(n,d,k,j)=
\\
\frac{\Nkjnums \probnums\, \EYdens} 
{ \Nkjdens \, \probdens \EYnums}
\end{multline*}
As usual, upon applying Stirling's formula we find that the powers of $n/e$
cancel and we are left with  
\begin{equation}\label{feq}
f(n,d,k,j)=f_0(n,d,k,j)G^n
\Bigpar{1+O\Bigpar{\frac1{j+1}+\frac1{k-j}+\frac1{n-k}}}
\end{equation}
 where  
 \begin{equation*}%\label{gdef}
f_0(n,d,k,j)= \frac{
 \sqrt{(dn-4n+2k)nd}}{2\pi(d-2)\sqrt{k (n-k)j(n-2k+j) }}
 \end{equation*}
 and
(with $\k= k/n$, $\gam=j/n$)
\begin{align}
G&= \frac{2^{\k-\gam}d^{d-2} (d-3)^{1+\gam-2\k}(d-4+2\k)^{d-4+2\k}
\k^\k(1-\k)^{3-3\k}} {(d-1)(d-2)^{
2d-6+\k+\gam}\gam^\gam (\k-\gam)^{2\k-2\gam}(1-2\k+\gam)^{1-2\k+\gam} }
\label{Gdef}
\\
&=F(\ga,\gd):=
\frac{2^\ga (t+2)^{t} 
(t-1)^{  \gd-\ga }(t-2\gd)^{t-2\gd} (1-\gd)^{1-\gd}\gd^{3\gd}} 
{(t+1)t^{
2t -2\gd  -\ga}(1-\gd-\ga)^{1-\gd-\ga}\ga^{2 \ga}
( \gd - \ga)^{\gd - \ga}}\nonumber
\end{align}
where $\gd=1-\k$, $\ga=\k-\gam$ and $t=d-2$. 

We seek the maximum value of $F$ in the triangle 
$$
T = \big\{(\gd,\ga):    0\le \ga\le \gd,\,  \gd+\ga\le 1\big\}.
$$
The partial derivatives of $\ln F$ are
$$
\frac{\partial (\ln F)}{\partial \gd} =
\ln\frac{ t^2    (t-1)\gd^3  (1-\gd-\ga) }{ (t-2\gd)^2(1-\gd)(\gd-\ga) },
$$
\begin{equation}\label{dela}
\frac{\partial (\ln F)}{\partial \ga} =
\ln\frac{2t(\gd-\ga)(1-\gd-\ga)      }{ (t-1)\ga^2 }.
\end{equation}
Setting these equal to zero gives necessary conditions for a
stationary point of $F$:
\begin{equation}\label{svante6}
t^2    (t-1)\gd^3  (1-\gd-\ga)=(t-2\gd)^2(1-\gd)(\gd-\ga)
\end{equation}
and
\begin{equation}\label{svante7}
(t+1)\ga^2-2t\ga+2t\gd(1-\gd)=0.
\end{equation}
Solving the first equation for $\alpha$ and substituting this into the
second shows that the value of $\gd$ at a stationary point 
in the interior of $T$  
must be a
root of     
\begin{equation}\label{gdef}
g(\gd) 
= (1-\gd)\big((t-2\gd)^2  -t^2(t-1)\gd^2\big)^2-2t^3\gd(1-2\gd)^2(t-2\gd)^2,
\end{equation}
which is quintic in $\delta$.
Now $g$  factorises as
\begin{equation}\label{gfact}
g(\delta)= (t+2)(\delta_0-\delta)h(\delta)
\end{equation}
where 
\begin{align*}
 h(\gd) =\  
 &t^5\gd^4-2\gd(2\gd^3+\gd^2-2\gd+1)t^4+(9\gd^4-12\gd^3+6\gd^2+1)t^3
\\
 &+2\gd(\gd-1)(3\gd^2+2\gd+3)t^2-4\gd^2(\gd+3)(\gd-1)t+8\gd^3(\gd-1),  
\end{align*}
and $\delta_0=t/(t+2)$ which, as we will show, determines the maximum
value of $F$ in $T$. 

The second derivative of $h$ with respect to $\delta$ is
\begin{align*}
&12t^5\delta^2+(-48\delta^2-12\delta+8)t^4+12(9\delta^2-6\delta+1)t^3\\
&+(72\delta^2-12\delta+4)t^2+(-48\delta^2-48\delta+24)t+96\delta^2-48\delta.
\end{align*}
This can be rearranged as 
\begin{align*}
&12(t-6)t^4\delta^2+(24\delta^2-12\delta+8)t^4+12(9\delta^2-6\delta+1)t^3
+12\delta(t^2-4) \\
&+ 12\delta^2(t^2-6t)+ 4(15\delta^2-6\delta+1)t^2+24(\delta^2-2\delta+1)t+96\delta^2,
\end{align*}
in which each collected term is clearly nonnegative for $t\ge 6$ and
all $\gd\ge0$. 
Thus $h$ is a convex function of $\gd$ for each $t$.  

 From~(\ref{gdef}) we compute firstly $g(1/2)=(t-1)^2(t-2)^4/32>0$, 
secondly 
$$
g(1/\sqrt t)<  t^2-2t^{5/2}(1-2/\sqrt t)^2(t-1)^2 
<0   
$$
since $t\ge 6$ and  $25(1-2/\sqrt 6)^2 >0.8$, and thirdly, since $t\ge
6$ implies $t^2-7t+8>0$ and thus $1-1/t<2(1-2/t)^2$, 
$$
g(1/t)< (1-1/t)\big((t-2/t)^4-t^2(t-2/t)^2\big)<0.
$$
Note that $h$ has the same sign as $g$ for $0\le \gd\le 1/2$ as
$\gd_0>1/2$. Thus $h(1/2)$ is positive and $h(1/\sqrt t)$ and $h(1/t)$
are negative. So by the convexity of $h$, the only zeros of $h$ for
$0<\gd < 1$  lie in $0< \gd <1/t$ or $1/\sqrt t < \gd < 1/2$. We
show separately that these two subsets of $T$ can hold no stationary
points of $F$.  
\smallskip

\noindent
{\bf Case 1: $0< \gd <1/t$}

\noindent
Substituting $\ga=\gb\gd$ into $\ln F$, and taking the second
derivative with respect to $\gd$, we obtain 
$$
\frac{\partial^2 \ln F(\gb\gd,\gd)}{\partial \gd^2} =
\frac{ 2t-\gb t-4t\gd +2t\gd^2-2\gb t\gd +2\gb t\gd ^2+2\gb \gd  }
   {( 1-\gd )( 1-\gd -\gb \gd )\gd ( t-2\gd )} 
$$
Since $\ga\le \gd$ in $T$, we have $\gb\le 1$ and so the factors in
the denominator are all positive. The numerator is at least 
$2t- t-4t\gd  -2  t\gd >0$  as $\gd<1/t\le 1/6$. Hence, $\ln F$  can
have no local maximum in $T$ for such $\gd$. 

\smallskip

\noindent
{\bf Case 2: $1/\sqrt t < \gd < 1/2$}

\noindent
For such $\gd$, from~(\ref{svante7}) we obtain $ 2t\ga>2t\gd(1-\gd)$,
i.e.\ $\ga> \gd(1-\gd)$. So, using~(\ref{svante6}), at a stationary
point 
$$
\frac{(1-  \gd)^2}{\gd^2}=
\frac{1- \gd- \gd(1-\gd)}{\gd-\gd(1-\gd)}< \frac{1-\gd -\ga}{\gd-\ga}
=\frac{(t-2\gd)^2(1-\gd) }{t^2(t-1)\gd^3   }
 $$
 and so, since $\gd^2>1/t$, 
 $$
 \frac{1-  \gd}{\gd}< \frac{(t-2\gd)^2  }{t(t-1)  }.
 $$
 But this fails at $\gd=1/2$, and the derivative of the left hand side
 with respect to $\gd$ is less than $-4$ for $\gd<1/2$, whilst that of
 the right is 
easily greater than $-4$. So the inequality fails, and there is no
 such stationary point. 
\smallskip

We conclude that $\gd=\gd_0$ determines the unique local maximum in
the interior of $T$. The boundary of $T$ must also be
investigated. Considering~(\ref{dela}),  
there is no local maximum at a boundary point with $0<\gd<1$, since
$\partial(\ln F)/\partial \ga$ tends to $\infty$ as $\ga$ tends to $0$
from the right, and to $-\infty$ as $\ga$ tends to $\min(\gd,1-\gd)$
from the left. A similar argument applies to eliminate
$(\ga,\gd)=(0,1)$ from consideration, since moving along the boundary
where $\ga=1-\gd$, $F$ is a smooth function times $\ga^{-\ga}$. This
leaves only the point $\ga=\gd=0$, which is indeed a local maximum,
with value $F(0,0)=(1+2/t)^t/(t+1)<1$ for $t\ge 6$.  

To deduce that $\gd_0$ determines the unique global maximum in $T$, we
only need to observe that the corresponding value of $\ga$ is
$\ga_0=2\gd_0/(t+1)$, and that $F(\ga_0,\gd_0)=1$. 

The  rest of the argument is totally standard for such variance
calculations, as in~\cite{FJMRW} for example, so we omit the
justifications. The point $(\ga_0,\gd_0)$ corresponds to $\kappa =
\kappa_0=2/d$, $\gam=\kappa_0/(d-1)$. Putting $\kappa=\kappa_0+\hat
\kappa/\sqrt n$ and $\gam=\gam_0+\hat \gam/\sqrt n$,   and expanding
$\ln(G^n)$ ($G$ defined in~(\ref{Gdef})) about $\hat \kappa=\hat \gam
= 0$, we find up to quadratic terms in $\hat \kappa$ and $\hat \gam$ 
$$
 G^n  \approx e^{c_1 \hat \kappa^2 + c_2 \hat \kappa\hat \gam +c_3\hat \gam^2}
$$
where 
\begin{align*}
c_1 &= -\frac{d(d^3-3d^2+4d+4)}{4(d-2)^2(d-3)}\\
c_2 & =   \frac { d (d-1)^2}{(d-2)(d-3)}\\
c_3 & =   -\frac{ d(d-1)^2}{4(d-3)}.
\end{align*}
Here $c_3$ is clearly negative, and the determinant $D=4c_1c_3-c_2^2$
 of the Hessian of the quadratic form is positive, as we expect since
 the expansion is at a local maximum. The routine argument now gives  
 from~(\ref{feq})
\begin{eqnarray*}
 \sum_{k,j<n} f(n,d,k,j)&\sim& f_0(n,d,n\kappa_0,n\gamma_0)2\pi n/\sqrt D \\
 &=&\sqrt\frac{d}{d-4}.
\end{eqnarray*}
Recalling~(\ref{mainvar}) and that $\E Y \to \infty$ now establishes
\refL{LE2}. 
\section{Rainbow matchings}

In this section we briefly consider the analoguous problem of the
existence of a rainbow perfect matching in a randomly coloured random
regular graph. We will omit the details of the calculations.

The model is now slightly different. We consider a random regular
graph $G(2n,d)$ with an even number of vertices ($d$ may now be arbitrary),
and colour randomly the $nd$ edges with $n$ colours, $d$ edges of each
colour.
We then ask whether there exists a rainbow matching consisting of $n$
disjoint edges of different colours.

We can translate this to a random bipartite (multi)graph as above; now
the bipartite graph has $2n$ plain vertices of degree $d$ and $n$
coloured vertices of degree $2d$.
Let $Z$ be the number of rainbow perfect matchings; in the bipartite version,
$Z$ is the number of decompositions of the graph into $n$ disjoint
paths of length 2, with 2 plain and 1 coloured vertex each.

Calculations as above yield
\begin{equation*}
  \E Z= d^{3n}\frac{(2n)!\,\bigpar{(2d-2)n}!}{(2dn)!}
=\Theta\xxpar{n\qq\parfrac{(d-1)^{2d-2}}{d^{2d-3}}^n}
\end{equation*}
and it is easily checked that, as \ntoo, 
$\E Z\to0$ for $d\le6$, while 
$\E Z\to\infty$ for $d\ge7$.
In particular, for $d\le6$ there is \whp{} no rainbow perfect
matching.

Furthermore, for $d\ge7$, an argument similar to the one in
\cite{BollMcK}
(and much simpler than the proof of \refL{LE2} above, since we only need
to maximize over one variable)
yields
\begin{equation*}
  \frac{\E(Z^2)}{(\E Z)^2} \to \frac{d-1}{\sqrt{d(d-3)}}.
\end{equation*}

Finally, defining $X\ij$ as before,  \refT{T0} applies when $d\ge7$
with
\begin{align*}
 \gl\ij&=\frac1{2i}\binom ij 2^j(d-1)^{i+j} ,
\\%\intertext{and}
 \gd\ij&=\frac{(-1)^{i+j}}{(d-1)^{i+j}}.
\end{align*}
Hence there exists a rainbow perfect matching \whp{} when $d\ge7$.

By analogy with the open problem in \refR{Rdecom} one might further  
ask whether there exists a
decomposition into $d$ rainbow 
perfect matchings.
Computing the expected number of such decompositions reveals that when
$d\le 11$, \whp{} no such decomposition exists. This compares with the
corresponding result for uncoloured graphs, that $G(2n,d)$ \whp{} has
a decomposition into $d$ perfect matchings as soon as
$d\ge3$~\cite{SJ103,MRRW}. 

\newcommand\AAP{\emph{Adv. Appl. Probab.} } %in?
\newcommand\JAP{\emph{J. Appl. Probab.} }
\newcommand\JAMS{\emph{J. \AMS} }
\newcommand\MAMS{\emph{Memoirs \AMS} }
\newcommand\PAMS{\emph{Proc. \AMS} }
\newcommand\TAMS{\emph{Trans. \AMS} }
\newcommand\AnnMS{\emph{Ann. Math. Statist.} }
\newcommand\AnnPr{\emph{Ann. Probab.} }
\newcommand\CPC{\emph{Combin. Probab. Comput.} }
\newcommand\JMAA{\emph{J. Math. Anal. Appl.} }
\newcommand\RSA{\emph{Random Strucures Algorithms} }
\newcommand\ZW{\emph{Z. Wahrsch. Verw. Gebiete} }
\newcommand\DMTCS{\jour{Discr. Math. Theor. Comput. Sci.} }

\newcommand\AMS{Amer. Math. Soc.}
\newcommand\Springer{Springer-Verlag}
\newcommand\Wiley{Wiley}

\newcommand\vol{\textbf}
\newcommand\jour{\emph}
\newcommand\book{\emph}
\newcommand\inbook{\emph}
\def\no#1#2,{\unskip#2, no. #1,} %(efter \rtal) eller

\newcommand\webcite[1]{\hfil\penalty0\texttt{\def~{\~{}}#1}\hfill\hfill}
\newcommand\webcitesvante{\webcite{http://www.math.uu.se/\~{}svante/papers/}}
\newcommand\arxiv[1]{arXiv:#1}

\def\nobibitem#1\par{}


\begin{thebibliography}{99}

\bibitem{bollobas}
B. Bollob{\'a}s, 
\book{Random Graphs}. 
Academic Press, New York, 1985; 
Second ed. Cambridge University Press, Cambridge, 2001. 

\bibitem{BollMcK}
B. Bollob{\'a}s \& B.D. McKay,
The number of matchings in random regular graphs and bipartite graphs.
\jour{J. Combin. Theory Ser. B} \vol{41} \no1 (1986),  80--91.


\bibitem{CF}
C. Cooper \& A. Frieze,
Multi-coloured Hamilton cycles in randomly coloured random graphs.
\CPC \vol{11}, 129--134.

\bibitem{CFMR}
C. Cooper, A. Frieze, M. Molloy \& B. Reed,
Perfect matchings in random $r$-regular, $s$-uniform hypergraphs.
\CPC \vol5 (1996), 1--15.

\nobibitem{friedman}
Friedman,~J. (1991) On the second eigenvalue and random walks in
random $d$-regular graphs, {\em Combinatorica} {\bf 11} 331--362.

\nobibitem{frieze}
Frieze,~A. (2001) Hamilton cycles in the union of random permutations,
{\em Random Structures Algorithms} {\bf 18} 83--94.

\bibitem{FJ}
A. Frieze \& S. Janson
Perfect matchings in random $s$-uniform hypergraphs.
\RSA \vol7 (1995), 41--57.

\bibitem{FJMRW}
A. Frieze, M. Jerrum, M. Molloy, R. Robinson \& N. Wormald,
Generating and counting Hamilton cycles in random regular graphs,
\emph{J. Algorithms} \vol{21} (1996), 176--198.

\bibitem{Gamburd}
A. Gamburd,
Poisson-Dirichlet distribution for random Belyi surfaces.
Preprint, 2005.
\arxiv{math.PR/0501283}

\bibitem{Garmo95}
H. Garmo, 
Random railways modeled as random $3$-regular graphs. 
\RSA \vol9 \no{1-2} (1996),  113--136. 

\bibitem{Garmo99}
H. Garmo, 
Asymptotic properties of the connectivity number of random railways. 
\AAP \vol{31} \no3 (1999), 720--741.
 		
\bibitem{GarmoJK}
H. Garmo, S. Janson \& M. Karo\'nski,
On generalized random railways.  
\CPC \vol{13}  \no1 (2004),  31--35.
	
\bibitem{GJKW}
C. Greenhill, S. Janson, J.H. Kim \& N.C. Wormald,
Permutation pseudographs and contiguity.
\CPC \vol{11} \no3 (2002),  273--298.

\bibitem{SJ103}
S. Janson,
Random regular graphs: asymptotic distributions and contiguity.
\CPC \vol4 (1995), 369--405.

\bibitem{JLR}
S. Janson, T. {\L}uczak \& A. Ruci{\' n}ski,
\book{Random Graphs}, \Wiley, New York,
2000.

\bibitem{KW99}
J.H. Kim \& N.C. Wormald,
Random matchings which induce Hamilton cycles,
and hamiltonian decompositions of random regular graphs,
\jour{J. Combin. Theory Ser. B} {\bf 81} (2001), 20--44.

\bibitem{MRRW}
M. Molloy, H. Robalewska, R.W. Robinson \& N.C. Wormald, 
1-factorisations of random regular graphs.
\RSA \vol{10} (1997), 305--321.

\nobibitem{rob}
Robalewska,~H. 2-factors in random regular graphs, {\em J. Graph
Theory} {\bf 23} 215--224.

\bibitem{RW92}
R.W. Robinson \& N.C. Wormald, 
Almost all cubic graphs are hamiltonian, 
\RSA \vol3 (1992), 117--125.

\bibitem{RW94}
R.W. Robinson \& N.C. Wormald,
Almost all regular graphs are hamiltonian, 
\RSA {\bf 5} (1994), 363--374.

\bibitem{RW01}
R.W. Robinson \& N.C. Wormald, 
Hamilton cycles containing randomly selected edges in random regular
graphs, 
\RSA {\bf 19} (2001), 128--147.

\bibitem{NW99}
N.C. Wormald,  
Models of random regular graphs.  
\inbook{Surveys in Combinatorics 1999},
eds. J.D. Lamb \& D.A. Preece,
LMS Lecture Note Series \vol{267},
Cambridge University Press, Cambridge, 1999,
239--298. 

\end{thebibliography}
\end{document}